\documentclass{article}
\usepackage{amsthm}
\usepackage{amssymb}
\usepackage{amsmath}
\usepackage[pdftex]{graphicx}
\usepackage{hyperref}
\hypersetup{pdftex,colorlinks=true,urlcolor=blue,pdfpagemode=none,pdfstartview=FitH,pdfview=FitH}
\usepackage{esdiff}
\usepackage{greek4cbc}
\usepackage{manfnt}
\usepackage{pifont}
\usepackage{mathrsfs}
\usepackage{pstricks}
\usepackage{ifsym}

\makeatletter
\def\revddots{\mathinner{\mkern1mu\raise\p@
      \vbox{\kern7\p@\hbox{.}}\mkern2mu
      \raise4\p@\hbox{.}\mkern2mu\raise7\p@\hbox{.}\mkern1mu}}
\makeatother

\begin{document}

\newtheorem*{Th}{Theorem}
\newtheorem*{Pro}{Proposition}
\newtheorem{De}[subsubsection]{D\'efinition}
\newtheorem{Prt}[subsection]{Propriété}
\newtheorem{Prts}[subsection]{Propriétés}
\newtheorem{Le}[subsubsection]{Lemme}
\newtheorem{Hyp}[subsection]{Hypothèse}

\newenvironment{desc}[1]{
	\begin{list}{}
		{\renewcommand{\makelabel}[1]{
			\textup{\textsf{##1}}\hfill}
		\settowidth{\labelwidth}{\textsf{#1 }}
		\leftmargin=\labelwidth
		\advance \leftmargin\labelsep
}}
	{\end{list}}

\bigskip

\title{A short proof of Cram\'er's theorem in $\mathbb{R}$}
\date{February 18, 2010}
\author{\\ Rapha\"el Cerf and Pierre Petit\\ \\
Universit\'e Paris Sud}
\maketitle

\bigskip

The most fundamental result in probability theory is the law of large numbers for a sequence $(X_n)_{n \geqslant 1}$ of independent and identically distributed real valued random variables. Define the empirical mean of $(X_n)_{n \geqslant 1}$ by
\[
\overline{X}_n = \frac{1}{n} \sum_{i=1}^n X_i
\]
The law of large numbers asserts that the empirical mean $\overline{X}_n$ converges (almost surely) towards the theoretical mean $\smash{\mathbb{E}( X_1 )}$ provided that 
$\smash{\mathbb{E}( |X_1| )}$ is finite.
The next fundamental results are the central limit theorem and Cram\'er's theorem. Both are refinements of the law of large numbers in two different directions. The central limit theorem describes the random fluctuations of $\overline{X}_n$ around $\smash{\mathbb{E}( X_1 )}$. Cram\'er's theorem estimates the probability that $\overline{X}_n$ deviates significantly from 
$\smash{\mathbb{E}( X_1 )}$ :
\[
\mathbb{P} \big( \overline{X}_n \geqslant \mathbb{E}( X_1 ) + \varepsilon \big)\textrm{\quad for $\varepsilon > 0$.}
\]
Such an event is called a "large deviations" event since it has a very small probability: it turns out that this probability decays exponentially fast with $n$. The first estimate of this kind can be traced back to Cram\'er's paper \cite{Cra38} which deals with variables possessing a density. In \cite{Che52} Chernoff relaxed this assumption. Then coming from statistical mechanics Lanford imported the subadditivity argument in the proof \cite{Lan73}. Cram\'er's theory was extended to infinite dimensional topological vector spaces by Bahadur and Zabell \cite{BZ79}. The classical texts of Azencott \cite{AZ}, 
Deuschel and Stroock \cite{DS89} and Dembo and Zeitouni \cite{DZ93} take stock of the foregoing improvements.
At this point in time classical proofs of Cram\'er's theorem in $\mathbb{R}$ resort either to the law of large numbers (see e.g. \cite{DZ93}) or to Mosco's theorem (see e.g. \cite{CRACER}). 
We expose here a direct proof of Cram\'er's theorem in $\mathbb{R}$ based on convex duality. Not only is the proof shorter, but it can easily adapt to a broader setting.


\begin{Th}[Cram\'er]
Let $(X_n)_{n \geqslant 1}$ be a sequence of independent and identically distributed real valued random variables and let $\overline{X}_n$ be the empirical mean:
\[
\overline{X}_n = \frac{1}{n} \sum_{i=1}^n X_i
\]
For all $x \in \mathbb{R}$ the sequence
\[
\frac{1}{n} \log \mathbb{P}\big( \overline{X}_n \geqslant x\big)
\]
converges in $[-\infty , 0]$ and
\[
\lim_{n\to\infty} \frac{1}{n} \log \mathbb{P}\big( \overline{X}_n \geqslant x\big) = \inf_{\lambda \geqslant 0} \big( \log \mathbb{E} \big( e^{\lambda X_1} \big) - \lambda x \big)
\]
\end{Th}


Let us define the entropy of the sequence $(X_n)_{n \geqslant 1}$ by
\[
\forall x \in \mathbb{R} \qquad s(x) = \sup_{n \geqslant 1} \frac{1}{n} \log \mathbb{P} \big( \overline{X}_n \geqslant x \big)
\]
and the pressure of $X_1$ (or the $\log$-Laplace transform of the law of $X_1$) by
\[
\forall \lambda \in \mathbb{R} \qquad p(\lambda) = \log \mathbb{E} \big( e^{\lambda X_1} \big)
\]
The entropy and the pressure may take infinite values.
Our strategy is to show a dual version, in the sense of convex functions, of Cram\'er's theorem.
\begin{Pro}[dual equality]
For all $\lambda \geqslant 0$,
\[
p(\lambda) = \sup_{u \in \mathbb{R}} \big( \lambda u + s(u) \big)
\]
\end{Pro}
\textbf{Proof.}
The classical Chebychev inequality will yield one part of the proof of the above equality. To prove the other part we condition $\overline{X}_n$ to be bounded by $K$ and then let $K$ grow towards $+ \infty$.
Since $X_1$, \ldots, $X_n$ are independent and identically distributed, we have for all $\lambda \geqslant 0$
\[
\forall n\geq 1\qquad
\mathbb{E} \big( e^{n \lambda \overline{X}_n} \big) = \mathbb{E} \big( e^{\lambda X_1} \big)^n
\]
Thus, for $u \in \mathbb{R}$ and $n\geq 1$, 
it follows from Chebychev inequality that
\begin{align*}
p(\lambda)  =  \log \mathbb{E} \big( e^{\lambda X_1} \big)
 & = \frac{1}{n} \log \mathbb{E} \big( e^{n \lambda \overline{X}_n} \big)\\
 & \geqslant_ç \frac{1}{n} \log \Big( e^{n\lambda u} \mathbb{P}\big( e^{n\lambda \overline{X}_n} \geqslant e^{n\lambda u} \big) \Big)
  \geqslant  \lambda u + \frac{1}{n} \log \mathbb{P} \big( \overline{X}_n \geqslant u \big)
\end{align*}
Hence, taking the supremum over $n \geqslant 1$ and then over $u \in \mathbb{R}$, we get
\[
\forall \lambda \geqslant 0 \qquad p(\lambda) \geqslant \sup_{u \in \mathbb{R}} \big( \lambda u + s(u) \big)
\]
Next, we prove the converse inequality for $\lambda = 0$ : for all $u \in \mathbb{R}$,
\begin{equation*}
0 \geqslant s(u) \geqslant \sup_{n \geqslant 1} \frac{1}{n} \log \mathbb{P} \big( X_1 \geqslant u \big)^n = \log \mathbb{P} \big( X_1 \geqslant u \big)
\end{equation*}
Hence, letting $u$ go to $-\infty$, we see that
\[
\sup_{u \in \mathbb{R}} s(u) = 0 = p(0)
\]
Now let $\lambda > 0$ and $K > 0$. For all $n \geqslant 1$, using the fact that $X_1$, \ldots, $X_n$ are independent and identically distributed, we have
\begin{align*}
\log \mathbb{E} \left( e^{\lambda X_1} 1_{|X_1| \leqslant K} \right) & = \frac{1}{n} \log \mathbb{E} \left( e^{\lambda (X_1 + \cdots + X_n)} 1_{|X_1| \leqslant K} \cdots 1_{|X_n| \leqslant K} \right) \\
 & \leqslant \frac{1}{n} \log \mathbb{E} \left( e^{n \lambda \overline{X}_n} 1_{|\overline{X}_n| \leqslant K} \right) \\
 & = \frac{1}{n} \log \mathbb{E} \left( \left( e^{-n\lambda K} + \int_{-K}^{\overline{X}_n} n\lambda e^{n \lambda u} du \right) 1_{|\overline{X}_n| \leqslant K} \right) \\
 & \leqslant \frac{1}{n} \log \left( e^{-n\lambda K} + \int \mathbb{E}\left( 1_{-K \leqslant u \leqslant \overline{X}_n} 1_{|\overline{X}_n| \leqslant K} \right) n\lambda e^{n\lambda u} du \right)
\end{align*}
the last step being a consequence of Fubini's theorem.
Since
\[
\mathbb{E}\left( 1_{-K \leqslant u \leqslant \overline{X}_n} 1_{|\overline{X}_n| \leqslant K} \right) \leqslant \mathbb{P} \big( \overline{X}_n \geqslant u \big) 1_{|u| \leqslant K} \leqslant e^{n s(u)} 1_{|u| \leqslant K}
\]
we get
\begin{align*}
\log \mathbb{E} \left( e^{\lambda X_1} 1_{|X_1| \leqslant K} \right) & \leqslant \frac{1}{n} \log \left( e^{-n\lambda K} + \int_{-K}^K n\lambda e^{n (\lambda u + s(u))} du \right)\\
 & \leqslant \frac{1}{n} \log \left( e^{-n\lambda K} + 2K n\lambda \exp \left( n \sup_{u \in \mathbb{R}} \big( \lambda u + s(u) \big) \right) \right)\\
\end{align*}
Let $K$ be large enough so that (recall that the supremum of $s$ is $0$)
\[
-\lambda K < \sup_{u \in \mathbb{R}} \big( \lambda u + s(u) \big)
\]
Sending $n$ to $\infty$ we obtain
\[
\log \mathbb{E} \left( e^{\lambda X_1} 1_{|X_1| \leqslant K} \right) \leqslant \sup_{u \in \mathbb{R}} \big( \lambda u + s(u) \big)
\]
Eventually sending $K$ to $+\infty$ we get
$$
\phantom{a}\qquad\qquad
\quad\qquad
\qquad\qquad
p(\lambda) \leqslant \sup_{u \in \mathbb{R}} \big( \lambda u + s(u) \big)
\quad\qquad
\qquad\qquad
\qquad\qquad
\qed
$$


To deduce Cram\'er's theorem from the dual equality we have just proved, we need some properties of the function $s$. First of all it follows from the definition of $s$ that $s$ is non-increasing. 
\begin{Pro}
For all $x \in \mathbb{R}$ the sequence
\[
\frac{1}{n} \log \mathbb{P}\big( \overline{X}_n \geqslant x\big)
\]
converges in $[-\infty , 0]$ towards $s(x)$. The function $s : \mathbb{R} \rightarrow [-\infty , 0]$ is concave.
\end{Pro}

\textbf{Proof.} Let $x, y \in \mathbb{R}$ with $x \leqslant y$. Let also $\alpha \in ]0, 1[$. Suppose that $\mathbb{P} \big( X_1 \geqslant y \big) > 0$ and let $n \geqslant m \geqslant 1$. Let $n = mq + r$ be the Euclidian division of $n$ by $m$. On the event
\[
\bigcap_{k=0}^{\lfloor \alpha q \rfloor-1} \Bigg\{ \sum_{i=mk+1}^{m(k+1)} X_i \geqslant mx \Bigg\} \cap \bigcap_{k=\lfloor \alpha q \rfloor}^{q -1}
\Bigg\{ \sum_{i=mk+1}^{m(k+1)} X_i \geqslant my \Bigg\} \cap \bigcap_{i=mq+1}^n \big\{ X_i \geqslant y \big\}
\]
we have (remember that $x \leqslant y$) :
\[
\sum_{i=1}^n X_i \geqslant \lfloor \alpha q \rfloor m x + \big( q - \lfloor \alpha q \rfloor \big) m y + r y \geqslant n \big( \alpha x + (1-\alpha) y \big)
\]
Therefore
\[
\mathbb{P}\big( \overline{X}_n \geqslant \alpha x + (1-\alpha) y \big) \geqslant \mathbb{P}\big( \overline{X}_m \geqslant x \big)^{\lfloor \alpha q \rfloor} \mathbb{P}\big( \overline{X}_m \geqslant y \big)^{q - \lfloor \alpha q \rfloor} \mathbb{P}\big( X_1 \geqslant y \big)^r
\]
Since $\mathbb{P}\big( X_1 \geqslant y \big) > 0$ and $0 \leqslant r < m$, taking logarithms, dividing by $n$ and sending $n$ to $\infty$, we get for all $m \geqslant 1$
\begin{align*}
\liminf_{n\to\infty} \frac{1}{n} \log \mathbb{P}\big( \overline{X}_n \geqslant \alpha x & + (1-\alpha) y \big)\\
& \geqslant \frac{\alpha}{m} \log \mathbb{P}\big( \overline{X}_m \geqslant x \big) + \frac{1-\alpha}{m} \log \mathbb{P}\big( \overline{X}_m \geqslant y \big)
\end{align*}
If $x = y$, taking the supremum over $m \geqslant 1$, we conclude that
\[
\lim_{n\to\infty} \frac{1}{n} \log \mathbb{P}\big( \overline{X}_n \geqslant x \big) = s(x)
\]
If $x < y$, sending $m$ to $\infty$, we get
\[
s \big( \alpha x + (1-\alpha) y \big) \geqslant \alpha s(x) + (1-\alpha) s(y)
\]
If $y$ is such that $\mathbb{P} \big( X_1 \geqslant y \big) = 0$, then
\[
\forall n \geqslant 1 \qquad \frac{1}{n} \log \mathbb{P}\big( \overline{X}_n \geqslant y\big) = - \infty
\]
whence $s(y) = - \infty$ and the concavity inequality still holds.
\qed

\medskip

We now finish the proof of Cram\'er's theorem.
At this point we know that for all $x \in \mathbb{R}$
\[
\inf_{\lambda \geqslant 0} \big( p(\lambda) - \lambda x \big) = \inf_{\lambda \geqslant 0} \sup_{u \in \mathbb{R}} \big( \lambda (u - x) + s(u) \big)
\]
It remains to prove that the latter quantity equals $s(x)$. The result is standard in convex functions theory. Let us give an elementary proof in our setting. The right-hand side of the previous equation is clearly superior or equal to $s(x)$ : take $u = x$. 
To prove the converse inequality 
we set
\[
c = \inf \big\{ x \in \mathbb{R} : \mathbb{P} \big( X_1 \geqslant x \big) = 0 \big\}
\]
and we distinguish the three cases $x < c$, $x > c$ and $x = c$.

$\bullet$ Suppose $x < c$. Since $s$ is concave and non-increasing let $- \lambda = s'_g(x) \leqslant 0$ be the left derivative of $s$ at point $x$. 
Then $s(x) > - \infty$ and
\[
\forall u \in \mathbb{R} \qquad s(u) \leqslant s(x) - \lambda(u - x)
\]
from which the result follows.

$\bullet$ Suppose $x > c$. Then $s(x) = -\infty$ and, for all $\lambda \geqslant 0$,
\[
p(\lambda) - \lambda x = \log \mathbb{E} \big( e^{\lambda(X_1 - x)} \big) \leqslant \log e^{\lambda (c - x)} = \lambda (c - x)
\]
so the infimum over $\lambda \geqslant 0$ is indeed $-\infty$.

$\bullet$ Suppose $x = c$. Then, for all $\lambda \geqslant 0$ and $\varepsilon > 0$,
\begin{align*}
p(\lambda) - \lambda c 
 & = \log \mathbb{E} \Big( e^{\lambda(X_1 - c)} \big( 1_{X_1 < c-\varepsilon} + 1_{c-\varepsilon \leqslant X_1 \leqslant c} \big) \Big)\\
 & \leqslant \log \Big( e^{-\lambda \varepsilon} + \mathbb{P} \big( X_1 \geqslant c - \varepsilon \big) \Big)
\end{align*}
Taking the infimum over $\lambda \geqslant 0$ and sending $\varepsilon$ to $0$ we get
$$
\phantom{a}\qquad\qquad\qquad\quad
\inf_{\lambda \geqslant 0} \big( p(\lambda) - \lambda x \big) \leqslant \log \mathbb{P} \big( X_1 \geqslant c \big) \leqslant s(c)
\quad\qquad\qquad\qquad\qed
$$

\vspace{-17pt}

\renewcommand\refname{}
\bibliographystyle{plain}
\bibliography{cramer}

\begin{thebibliography}{1}

\bibitem{AZ}
R.~Azencott.
\newblock {\em Grandes d\'eviations et applications}.
\newblock Ecole d'Et\'e de Probabilit\'es de Saint--Flour, Lecture Notes in
  Mathematics 774. Springer--Verlag, 1980.

\bibitem{BZ79}
R.~R. Bahadur and S.~L. Zabell.
\newblock Large deviations of the sample mean in general vector spaces.
\newblock {\em Ann. Prob.}, 7(4):587--621, 1979.

\bibitem{CRACER}
R.~Cerf.
\newblock {\em On {C}ram\'er's theory in infinite dimensions}.
\newblock Panoramas et Synth\`eses 23. Soci\'et\'e Math\'ematique de France,
  Paris, 2007.

\bibitem{Che52}
H.~Chernoff.
\newblock A measure of asymptotic efficiency for tests of a hypothesis based on
  the sum of observations.
\newblock {\em Ann. Math. Statist.}, 23:493--507, 1952.

\bibitem{Cra38}
H.~Cram\'er.
\newblock Sur un nouveau th\'eor\`eme limite de la th\'eorie des
  probabilit\'es.
\newblock {\em Actualit\'es scientifiques et industrielles}, 736:5--23, 1938.

\bibitem{DZ93}
A.~Dembo and O.~Zeitouni.
\newblock {\em Large deviations techniques and applications}.
\newblock Springer--Verlag, 1998.

\bibitem{DS89}
J.-D. Deuschel and D.~W. Stroock.
\newblock {\em Large deviations}.
\newblock Academic Press, 1989.

\bibitem{Lan73}
O.~E. Lanford.
\newblock {\em Entropy and equilibrium states in classical statistical
  mechanics}.
\newblock Lecture Notes in Physics 20. Springer, 1973.

\end{thebibliography}

\end{document}